\documentclass[preprint,12pt]{elsarticle}

\usepackage{amssymb}
\usepackage{amsmath}
\usepackage{amsthm}
\usepackage{vmargin}
\setmarginsrb{3.2cm}{1cm}{3.2cm}{3cm}{1cm}{1cm}{1.5cm}{1.5cm}
\newtheorem{thm}{Theorem}
\newtheorem{prop}{Proposition}
\newtheorem{lem}[thm]{Lemma}
\newdefinition{rem}{Remark}
\newdefinition{defi}{Definition}
\newproof{pf}{Proof}
\newproof{pot}{Proof of Theorem \ref{thm2}}
\def\arccos{\operatorname{arccos}}

\journal{*****}

\begin{document}

\begin{frontmatter}

\title{Non-integrability of a self-gravitating Riemann liquid ellipsoid}
\author{Thierry COMBOT\fnref{label2}}
\ead{combot@imcce.fr}
\address{IMCCE, 77 Avenue Denfert Rochereau 75014 PARIS}

\title{Non-integrability of a self-gravitating Riemann liquid ellipsoid}

\author{}

\address{}

\begin{abstract}
We consider the motion of a triaxial Riemann ellipsoid of homogeneous liquid without angular momentum. We prove that it does not admit an additional first integral which is meromorphic in position, impulsions, and the elliptic integrals which appear in the potential. This proves that the system is not integrable in the Liouville sense; we actually show that even its restriction to a fixed energy hypersurface is not integrable.
\end{abstract}

\end{frontmatter}

\section{Introduction}

We consider the motion of a self gravitating triaxial Riemann ellipsoid of homogeneous liquid, in the restricted case of zero angular momentum. Noting
$$J(q_1,q_2)=\int\limits_{0}^\infty \frac{1}{\sqrt{\left(\eta+\frac{4}{q_2^2}\right)\left(\eta^2+r\eta+\frac{q_2^2}{4}\right)}} d\eta$$
with $r^2=q_1^2+q_2^2$, this problem has the following Hamiltonian formulation
$$H=r\left(p_1^2+\frac{q_2^4p_2^2}{q_2^4+r}\right)+\alpha J(q_1,q_2) $$
This equation can be found in \cite{28}. The motion of Riemann ellipsoids have been studied for a long time, as for example in \cite{29},\cite{30}. {The integral $J$ and the Hamiltonian $H$ are regarded from a dynamical point of view as real-valued functions of real variables. In the case of $\alpha=0$, Ziglin manage in \cite{31} to extend the real system to the complex domain, which allowed him to prove the meromorphic non-integrability of the Hamiltonian $H$}. This case can be seen as the limit case when we take an energy level which tends to infinity. Still, a first integral of $H$ could be non-meromorphic when energy tend to infinity. Thus, for $\alpha\neq 0$, this does not prove ``meromorphic non-integrability'' of $H$ in a reasonable sense, nor does it give a reasonable notion of ``meromorphic non-integrability'' for $H$, as $H$ is not meromorphic in $p,q$ itself. {Indeed, the integral $J$ is multivalued in the complex domain. Thus, when $\alpha\neq 0$, extending $H$ in the complex requires to build a covering $\mathcal{M}$ on which the Hamiltonian is single valued. Such work is necessary to apply then the Morales-Ramis Theorem \cite{11}}.

As the case $\alpha=0$ has been already treated in \cite{31}, we will only consider the case $\alpha\neq 0$, which, after variable change, reduces to $\alpha=1$. The Hamiltonian $H$ has two degrees of freedom, and there is an elliptic integral in the potential. Because of this, we need to be precise about the notion of integrability we want to study. The main Theorem of this article is the following
\begin{thm}\label{thm0}
The Hamiltonian $H$ does not admit a first integral which is meromorphic on $\mathcal{M}$ and functionally independent with the Hamiltonian $H$. On a fixed energy hypersurface $H=h$, the Hamiltonian system restricted to this hypersurface does not admit a non-constant first integral meromorphic on $\mathcal{M}$.\\
\end{thm}
The main difficulty is that, contrary to algebraic functions, the transcendental extension $J$ cannot be suppressed in some additional first integral only by algebraic transformations. Instead, we will use a similar approach done by this author in \cite{91} to deal with algebraic extensions in potentials. The structure of the proof is the following
\begin{itemize}
\item First we build the manifold $\mathcal{M}$ on which $H$ is well defined.
\item We compute a family of particular solutions corresponding to the case when the ellipsoid of fluid is invariant by rotation. Closed form solutions can be given in this case, and we analyze (complex) singularities of these solutions.
\item We compute the variational equation near these orbits and prove a small extension of Morales-Ramis for proving the non existence of an additional first integral which is meromorphic on $\mathcal{M}$.
\item As the variational equation has transcendental functions in its coefficients, Kovacic algorithm \cite{78} cannot be used directly. Instead, we prove using an analysis of the monodromy that some second order differential equation related to the variational equation should have a Galois group whose identity component is solvable (which is a weaker property than the classical virtual abelianity condition from Morales-Ramis).
\end{itemize}
By chance, even if the Hamiltonian has elliptic functions, in the computation we will never see anything ``worse'' than that. The difficulty is that computing Galois group with transcendental functions in the base field can be much more tricky than in the case of rational base field.

\section{The phase manifold}

We may compute a closed form for the integral $J$ appearing in the Hamiltonian. Let $w=(q_2^6-16rq_2^2+64)^{1/4}$. Then
\begin{equation}\begin{split}\label{eq000}
J=\frac{\sqrt{2}\, q_2}{w} \left(2\mathcal{K}\left(\sqrt{\frac{w^2-rq_2^2+8}{2w^2}}\,\right)-
\mathcal{F}\left(\frac{4\sqrt{2}w}{8+w^2},\sqrt{\frac{w^2-rq_2^2+8}{2w^2}}\,\right)\right)\\
\end{split}\end{equation}
where $\mathcal{K}$ and $\mathcal{F}$ correspond respectively to the complete and incomplete elliptic integral of the first kind. Moreover we have
$$sn(2\mathcal{K}(\sqrt{k}\,)-(2\mathcal{K}(\sqrt{k}\,)-\mathcal{F}(z,\sqrt{k}\,)),\sqrt{k}\,)=z$$
where $sn$ denotes the Jacobi elliptic function $sn$. The function $sn$ is meromorphic in its two variables and the function $\mathcal{K}(\sqrt{k}\,)$ is well defined and holomorphic on the unit open disk $D(0,1)$. Let us now build a complex manifold on which the Hamiltonian $H$ will be well defined.

\begin{defi}
We consider the following functions on $\mathbb{C}^5\times D(0,1)$
\begin{equation}\begin{split}\label{eqrela}
g_1(q_1,q_2,r,w,k,j) & =r^2-q_1^2-q_2^2\\
g_2(q_1,q_2,r,w,k,j) & =q_2^6-16rq_2^2+64-w^4 \\
g_3(q_1,q_2,r,w,k,j) & =w^2-rq_2^2+8-2kw^2 \\
g_4(q_1,q_2,r,w,k,j) & =(w^2+8) sn(2\mathcal{K}(\sqrt{k}\,)-j,\sqrt{k}\,)-4\sqrt{2}\, w 
\end{split}\end{equation}
Let $E\subset \mathbb{C}^7\times D(0,1)$ be the set on which the Jacobian matrix of the application $(r,w,j,k)\longrightarrow g(r,w,k,j)$ is invertible. We define the \emph{phase manifold} to be
\begin{equation*}\begin{split}
\mathcal{M}=\left\lbrace (p,q,r,w,j,k)\in E,\;\; g_i(q_1,q_2,r,w,j,k)=0,\;i=1\dots 4, \; w(q_2^4+r)\neq 0 \right\rbrace
\end{split}\end{equation*}
\end{defi}

\begin{prop}\label{thm01}
The Hamiltonian $H$ is a holomorphic function (univalued) on $\mathcal{M}$, which is a complex analytic symplectic manifold.
\end{prop}

\begin{proof}
In the definition of $\mathcal{M}$, we look at the constraints on $(q_1,q_2,r,w,j,k)$. The functions $g_i$ are meromorphic (and univalued) on $\mathbb{C}^5\times D(0,1)$. These constraints are functionally independent (because in each constraint there is a variable that does not appear in the previous ones), so the rank of the associated Jacobian matrix is maximal outside possibly the set $E$ which is at least of codimension $1$. So $\mathcal{M}$ is a complex analytic manifold of dimension $4$, and we put on $\mathcal{M}$ the canonical symplectic structure in $p,q$ (which is well defined when $(p,q,r,w,j,k)\in E$, as proved by this author in \cite{91}). The integral $J$ and the Hamiltonian $H$ are well defined (univalued) on $\mathcal{M}$ as they can be expressed, for $(p_1,p_2,q_1,q_2,r,w,k,j)\in\mathcal{M}$, as
$$J=\frac{\sqrt{2}\,q_2 j}{w}\qquad H=r\left(p_1^2+\frac{q_2^4p_2^2}{q_2^4+r}\right)+\frac{\sqrt{2}\,q_2 j}{w}$$
These expressions are rational. The singularities of $H$ correspond to $w(q_2^4+r)=0$, which is outside of $\mathcal{M}$. Thus the Hamiltonian $H$ is a holomorphic function on $\mathcal{M}$.
\end{proof}

\begin{defi} (see Shafarevich \cite{27} p 362)
A representant of meromorphic function on $\mathcal{M}$ is given by an atlas of connected open sets $U_{\alpha,\;\alpha\in I}$ covering the manifold $\mathcal{M}$ and couples of holomorphic functions $(h_\alpha,k_\alpha)$ on each open set $U_\alpha$ such that
\begin{itemize}
\item The functions $k_\alpha$ are not identically $0$ on $U_\alpha$.
\item For all $\alpha,\beta \in I$, we have $h_\alpha k_\beta=h_\beta k_\alpha$ on $U_\alpha\cap U_\beta$.
\end{itemize}
We say that two representants $f_1,f_2$ of meromorphic functions on $\mathcal{M}$
$$f_1=\{(U_\alpha,h_\alpha,k_\alpha),\;\alpha\in I_1\} \qquad f_1=\{(\tilde{U}_\beta,\tilde{h}_\beta,\tilde{k}_\beta),\; \beta\in I_2\}$$  
are equivalent if for all $(\alpha,\beta)\in I_1\times I_2$, we have $h_\alpha \tilde{k}_\beta=\tilde{h}_\beta k_\alpha$ on $U_\alpha \cap \tilde{U}_\beta$.
Finally, a meromorphic function on $\mathcal{M}$ is a representant of meromorphic function on $\mathcal{M}$ quotiented by this equivalence relation.
\end{defi}

The meromorphic functions on $\mathcal{M}$ form now a field, and will be noted $\mathcal{M}er(\mathcal{M})$. The process of localization in the definition of meromorphic functions is necessary as the ring of holomorphic functions on an analytic manifold is not an integral domain, as in the following example
$$f_1,f_2:\Delta=\{(x,y)\in\mathbb{C}^2, xy=0,\; (x,y)\neq (0,0)\} \longrightarrow \mathbb{C}$$
$$f_1(x,y)=x\quad f_2(x,y)=y$$
The product $f_1f_2$ is identically $0$ on the analytic manifold $\Delta$.

\begin{defi}\label{def2}
We set the following definitions
\begin{itemize}
\item An element $I\in \mathcal{M}er(\mathcal{M})$ is a \emph{meromorphic first integral} of $H$ if $I$ is constant along any orbit.
\item An element $I\in \mathcal{M}er(\mathcal{M})$ is a \emph{meromorphic first integral on the level $H=h$} if $I$ is constant along any orbit with energy $H=h$.
\item An element $I\in \mathcal{M}er(\mathcal{M})$ is an \emph{additional first integral of $H$} (respectively on the level $H=h$), if $I$ is functionally independent with $H$ (respectively non-constant on $H=h$).
\end{itemize}
\end{defi}

\begin{rem}
The field of meromorphic functions $\mathcal{M}er(\mathcal{M})$ is a differential field (with respect to the derivations $\partial_{q_1},\partial_{q_2},\partial_{p_1},\partial_{p_2}$), and contains the Hamiltonian. We should also notice that even if a first integral on a level $H=h$ is defined a priori everywhere according to Definition \ref{def2}, in fact we only care about its restriction on $H=h$.
\end{rem}

The ``tool'' we will use to prove non-integrability is the Morales-Ramis theory from \cite{20},\cite{11}, and its extensions like \cite{91}.

\begin{thm}\label{thmmorales} (Morales-Ramis  \cite{11})
Let $H$ be a Hamiltonian holomorphic on a complex symplectic manifold $M$ of dimension $2n$, and $\Gamma \subset M$ a non-stationary orbit of $H$. If there are $n$ meromorphic first integrals of $H$ that are in involution and independent over a neighbourhood of $\Gamma$, then the identity component of the Galois group of the variational equation near $\Gamma$ is abelian.
\end{thm}

If the Hamiltonian satisfies good properties (as being holomorphic along the curve $\Gamma$), then a meromorphic first integral of $H$ will produce a meromorphic initial form on $\Gamma$ and then a constraint on the Galois group. This is a very general way to prove non-integrability: If a Hamiltonian is integrable in some sense, then the Galois group of the variational equation near a particular orbit should satisfy some particular property. In their article \cite{11}, they consider meromorphic Hamiltonians and first integrals, and here we will also prove an extension of Morales-Ramis Theorem to include elliptic integrals.\\

\section{Variational equation}

\subsection{An invariant manifold}

To study integrability, we will need to look at a particular (explicit) solution. The solution we will study corresponds in fact to the case where the ellipsoid of fluid is invariant by rotation.

\begin{prop}\label{thm1}
We consider the following set
\begin{equation}\begin{split}
\mathcal{P}=\{(p,q,r,w,k,j)\in \mathcal{M},
(p_1,q_1,k,r-q_2,q_2^3(\cos j-1)+16)=0 \}
\end{split}\end{equation}
The Hamiltonian vector field $X_H$ on $\mathcal{P}$ is tangent to $\mathcal{P}$.
\end{prop}

\begin{proof}
The manifold $\mathcal{M}$ and the Hamiltonian $H$ are invariant under $q_1\longrightarrow -q_1$. So $q_1=0$ is a plane of symmetry and thus is invariant. We now fix $q_1=0$ in the relations \eqref{eqrela} and we get the constraints
$$r^2-q_2^2 =0 \qquad  q_2^6-16rq_2^2+64-w^4 =0 \qquad w^2-rq_2^2+8-2kw^2 =0 $$
$$(w^2+8) sn(\mathcal{K}(\sqrt{k}\,)-j,\sqrt{k}\,)-4\sqrt{2}\, w =0$$
Let us now note
$$sn(\mathcal{K}(\sqrt{k}\,)-j,\sqrt{k}\,)=s,\;\; cn(\mathcal{K}(\sqrt{k}\,)-j,\sqrt{k}\,)=-c$$
This produces a polynomial ideal $I=$
$$<\! r^2-q_2^2,q_2^6-16rq_2^2+64-w^4,w^2-rq_2^2+8-2kw^2,(w^2+8) s-4\sqrt{2}\, w,c^2+s^2-1\!>$$
{This ideal is not prime. Let us note $\mathcal{Z}_1,\dots,\mathcal{Z}_n$ the algebraic varieties associated to each prime component of this ideal. We have the following inclusion
$$(\mathcal{M} \cap \{q_1=0\}) \subset \bigcup\limits_{i=1}^n \mathcal{Z}_i$$
The Hamiltonian field $X_H$ on $\mathcal{M} \cap \{q_1=0\}$ is a holomorphic vector field in the tangent space of the algebraic manifold $\cup_{i=1}^n \mathcal{Z}_i$, and thus for each $i=1\dots n$, we have}
$$\left. X_H\right|_{(\mathcal{M} \cap \{q_1=0\}\cap \mathcal{Z}_i)} \subset \mathcal{TZ}_i$$
Performing a prime decomposition of the ideal $I$ (using the command \textbf{PrimeDecomposition} of Maple), we then consider the prime factor of ideal $I$
\begin{equation}\begin{split}\label{eqrelP}
<\! r-q_2,k,q_2^3 s-4\sqrt{2}\,w,c^2+s^2-1,q_2^3c-q_2^3+16\!>
\end{split}\end{equation}
{and we note $\mathcal{Z}_1$ the associated algebraic variety. On $\mathcal{M}\cap \{q_1=0\}\cap \mathcal{Z}_1$ the Hamiltonian vector field is tangent to $\mathcal{M}\cap \{q_1=0\}\cap \mathcal{Z}_1$. We have moreover $k=0$, and then the elliptic functions $sn,cn$ become $s=\sin j, c=\cos j$. Thus $\mathcal{P}=\mathcal{M}\cap \{q_1=0\}\cap \mathcal{Z}_1$}.
\end{proof}

\begin{prop}\label{thm12}
The set $\mathcal{P}$ is a $2$-dimensional analytic manifold.
\end{prop}

\begin{proof}
{The set $\mathcal{P}\subset \mathbb{C}^5\times D(0,1)$ is given by the following equations
$$\mathcal{P}=\{(p,q,r,w,k,j)\in E,\;\;w(q_2^4+r)\neq 0,\;$$
$$(p_1,q_1,k,r-q_2,q_2^3(\cos j-1)+16,q_2^3\sin j -4\sqrt{2}\, w)=0\}$$
as the other equations in the definition of $\mathcal{M}$ are implied by these ones. The right equations are closed conditions, and the conditions $(p,q,r,w,k,j)\in E$, $w(q_2^4+r)\neq 0$ are open conditions. Let us first check that $\mathcal{P}$ is an analytic manifold. We compute the Jacobian matrix (for the $3$ last constraints in $q_2,r,j,w$)
$$J=\left( \begin {array}{cccc} 
-1&1&0&0\\3q_2^2(\cos j-1)&0&-q_2^3\sin j&0\\3q_2^2\sin j&0 &q_2^3\cos j &-4\sqrt{2}\\
\end {array} \right) $$
Let us look at the rank of this matrix. Taking the determinant of the three last columns, which is
$4\sqrt{2}\,q_2^3\sin j$, we obtain that the rank of $J$ is $3$ except maybe for $w=0$. But $w=0$ is not allowed in $\mathcal{P}$. Thus the rank of the Jacobian matrix $J$ of the closed conditions is always $3$, and thus $\mathcal{P}$ is an analytic manifold (recall that the closed conditions are analytic).}

Let us now look at the set $E$. The dimension of $\mathcal{P}$ could be lower than expected if the open conditions exclude most of the points satisfying the closed conditions. The set $E$ is the set of points $(p,q,r,w,j,k)$ such that the Jacobian matrix of the application $(r,w,j,k)\longrightarrow g(q_1,q_2,r,w,j,k)$ is invertible. So we need to study the rank of this Jacobian matrix under the conditions
\begin{equation}\label{eq01}
(p_1,q_1,k,r-q_2,q_2^3(\cos j-1)+16,q_2^3\sin j -4\sqrt{2}\, w)=0
\end{equation}
The Jacobian matrix of the application $(r,w,j,k)\longrightarrow g(q_1,q_2,r,w,j,k)$ under the conditions \eqref{eq01} are given by
$$\left( \begin {array}{cccc} 
2r&0&0&0\\-16q_2^2&-4w^3&0&0\\-q_2^2&2w&-2w^2&0\\0&2w\sin j-4\sqrt{2}&f(j)&(w^2+8)\cos j\\
\end {array} \right) $$
The determinant of this matrix is
$$\hbox{det}=8\cos j wq_2(w^2+8)(q_2^2w^2+2w^2r^2-3q_2^4r+8q_2^2+16r^2)$$
Noting $c=\cos j$, we build the ideal describing the singular locus, generated by equations \eqref{eq01}, $\hbox{det}\,w\,(q_2^4+r)$, $c^2+s^2-1$. This ideal is zero-dimensional. {Thus there are at most a discrete set of points satisfying the closed conditions \eqref{eq01} which are not in $E$ or such that $w(q_2^4+r)= 0$. As these open conditions remove finitely many points, the dimension of the manifold given by equations \eqref{eq01} is not reduced, and thus $\mathcal{P}$ is a two dimensional analytic manifold.}
\end{proof}

\begin{prop}
The restriction of the Hamiltonian $H$ to $\mathcal{P}$ gives the Hamiltonian
$$R(p_2,q_2)=\frac{q_2^4p_2^2}{q_2^3+1}+ \frac{\sqrt{2}\,q_2 j}{w} $$
On a each level $R=h$, there are finitely many critical points and a non-stationary orbit.
\end{prop}

\begin{proof}
The function $R$ is easily found by direct computation. The function $R$ defines a one degree of freedom Hamiltonian on $\mathcal{P}$. So, the orbit can be completely studied analyzing the levels of $R$.

Let us first prove that $R=h$ is a curve (not simply points). We consider generators of \eqref{eqrelP} and the equation $R-h=0$. On each equation, there is a variable which does not appear in the previous ones. Thus these constraints are functionally independant, and so $R=h$ is a curve.

We now look for critical points. For a critical point, we have
$$\dot{q}_2=\frac{2q_2^4p_2}{q_2^3+1}$$
As in the ideal \eqref{eqrelP} defining $\mathcal{P}$, $q_2$ cannot vanish, this implies that a critical point always corresponds to $p_2=0$. Differentiating the generators of \eqref{eqrelP} on $\mathcal{P}$, we get the relations
$$-4\sqrt{2}\,w\frac{\partial j}{\partial q_2}+3q_2^2c-3q_2^2=0 \qquad  2w\frac{\partial w}{\partial q_2} -3q_2^2=0$$
Now we differentiate $R$ in $q_2$ for $p_2=0$ and this produces
\begin{equation}\begin{split}\label{eq1}
\left. \frac{\partial R}{\partial q_2}\right|_{p_2=0} & = \frac{\sqrt{2}\, j}{w}+ \frac{\sqrt{2}\, q_2}{w}\frac{\partial j}{\partial q_2}-\frac{\sqrt{2}\,q_2 j}{w^2}\frac{\partial w}{\partial q_2}\\
& =\frac{\sqrt{2}\, j}{w} -\frac{3\sqrt{2}\,q_2^3 j}{2w^3}+\frac{3q_2^3(c-1)}{4w^2}\\
& =-\frac{j(q_2^3+16)}{\sqrt{2}\,w(q_2^3-8)}+\frac{3q_2^3(c-1)}{4w^2}
\end{split}\end{equation}
We need now to solve the equation $(\partial_{q_2} R,R-h)=0$ on $\mathcal{P}$. This corresponds to add to ideal \eqref{eqrelP} the ideal
$$<\!\sqrt{2}\,q_2j-hw, 4\sqrt{2}\,jw^2-6\sqrt{2}\,jq_2^3+3q_2^3wc-3q_2^3w\!>$$
Eliminating in $h,j$, we obtain the condition $27h^6j+4j^7+36j^5-216h^3j=0$. So for each $h$, there are finitely many possible $j$, and thus finitely many possible $r,q_2,w,c,s$.
\end{proof}

\subsection{Normal variational equation}

\begin{lem}\label{thm2}
The normal variational equation near a non-stationary orbit in $\mathcal{P}$ with $R=h$ is given by
\begin{equation}\begin{split}\label{eq2}
16(h-q_2\delta)(q_2^3+1)(q_2^3-8)^2q_2^4X''+\\
4(q_2^7\delta +5q_2^4\delta +112q_2\delta +24q_2^4+24q_2+12hq_2^3-96h)(q_2^3-8)q_2^3X'-\\
(q_2^3+1)^2(q_2^6\delta -8q_2^3-32q_2^3\delta -128)X=0
\end{split}\end{equation}
with $'$ corresponding to the derivation in $q_2$ and $\delta =\sqrt{2}\,jw^{-1}$.
\end{lem}

\begin{proof}
We begin by direct computations of the Hessian matrix of $H$ on $\mathcal{P}$. This works easily in all cases except for $\partial_{q_1q_1} H,\partial_{q_2q_2} H$. We obtain in particular that the variational equation, with $Y=(\Delta p_1,\Delta q_1,\Delta p_2,\Delta q_2)$, is of the form
$$\dot{Y}=\left(\begin{array}{cccc}0&-F_1&0&0\\2q_2&0&0&0\\ 0&0&*&*\\ 0&0&*&*\\ \end{array}\right) Y$$
So the variational equation is already decoupled, and the normal part corresponds to the first $2\times 2$ block. We have written here $F_1=\partial_{q_1q_1} H$. Let us look now closer at this function. The kinetic part after differentiation disappears on $\mathcal{P}$. So the only thing to compute is $\partial_{q_1q_1} \left(\sqrt{2}\, q_2 j/w\right)$ on $\mathcal{P}$.

We first need to compute differential relations in $q_1$ of $r,w,k,j$ on $\mathcal{P}$. We use the formula
$$\left.\partial_{kk} sn(z,k)\right|_{k=0}=\textstyle{\frac{1}{2}} \cos z (\sin z \cos z-z)$$
We differentiate the constraints \eqref{eqrela} and evaluate them on $\mathcal{P}$.
$$\partial_{q_1} r=0\quad \partial_{q_1} w=0 \quad \partial_{q_1} j=0 $$
At second order, we get (remember the relations \eqref{eqrelP})
$$\left.\partial_{q_1q_1 }r\right|_{q_1=0}=q_2^{-1}\quad \left.\partial_{q_1q_1 }w\right|_{q_1=0} = -4q_2 w^{-3} \quad  \left.\partial_{q_1q_1 }k\right|_{q_1=0} = -\frac{q_2^4}{2(q_2^3-8)^2}$$
$$\left.\partial_{q_1q_1 }j\right|_{q_1=0}=
-\frac{q_2^4j}{8w^4}+\frac{(q_2^3+16)\sqrt{2}}{2q_2^2w^3}$$
Using these formulas, we get the expression
$$F_1=-\frac{\sqrt{2}q_2^2(q_2^3-32)j}{8w^5}+\frac{q_2^3+16}{q_2w^4} $$
We now transform the normal variational equation in a second order differential equation. This gives (the dot being the derivation in time)
$$q_2(t)\ddot{X}-\dot{q_2}(t)\dot{X} = -2q_2(t)^2F_1(q_2(t))X$$
for an orbit of $H$ in $\mathcal{P}$. We also get that
\begin{equation}\label{eqphi1}
p_2(t)=\frac{\dot{q_2}(t)(q_2(t)^3+1)}{2q_2(t)^4}
\end{equation}
We now want to make the variable change $q_2(t)\longrightarrow t$. For this, we use the fact that $R$ is constant over an orbit and we get the relation
\begin{equation}\label{eqphi2}
\dot{q_2}(t)^2=\frac{4q_2(t)^4\left(hw-\sqrt{2}\,q_2(t) j\right)}{w(q_2(t)^3+1)}
\end{equation}
and then a relation for $\ddot{q_2}(t)$. The variable change $q_2(t)\longrightarrow q_2$ involves only $\ddot{q_2}(t)$ and $\dot{q_2}(t)^2$. This eventually produces equation \eqref{eq2}.
\end{proof}

\begin{rem}
The extension $\delta$ appearing in the normal variational equation can be written as a (multivalued) function of $q_2$
$$\delta =\frac{i \arccos\left(2\sqrt{2}q_2^{-3/2}\right)}{\sqrt{2} \sqrt{8-q_2^3}}$$
This variational equation does not have rational coefficients, and is in fact well defined on a Riemann surface.
\end{rem}

\section{Non-integrability}

Let $\mathcal{M}er(\mathcal{S})$ be the field of meromorphic functions on the Riemann surface
$$\mathcal{S}=\{(q_2,j)\in\mathbb{C}^2,\;\;q_2^3(\cos j-1)+16=0,\;\;(q_2^3+1)(q_2^3-8)\neq 0 \}$$
This differential field will be the base field for Galois group computations.

\begin{thm}\label{thm3}
If $H$ has an additional first integral in $\mathcal{M}er(\mathcal{M})$, then the Galois group over the base (differential) field $\mathcal{M}er(\mathcal{S})$ of the normal variational equation \eqref{eq2} is virtually Abelian for any $h\in\mathbb{C}$. If $H$ has an additional first integral in $\mathcal{M}er(\mathcal{M})$ on some fixed energy level $H=h$, then the Galois group over $\mathcal{M}er(\mathcal{S})$ of the normal variational equation \eqref{eq2} is virtually Abelian.
\end{thm}

This theorem is in fact a small extension of the Morales-Ramis Theorem \ref{thmmorales} and of the Ayoul-Zung Theorem \cite{81}; it is very similar to the one given by this author in \cite{91} for algebraic functions.

\begin{proof}
First remark that the coefficients of equation \eqref{eq2} are in $\mathcal{M}er(\mathcal{S})$ because we have
$$s=\frac{8j}{q_2^3\sin j}\;\;\in\mathcal{M}er(\mathcal{S})$$
The Hamiltonian $H$ is defined on the complex symplectic manifold $\mathcal{M}$. We consider the curve given by
$$\Gamma_h=\mathcal{P}\cap \{R=h\} \subset \mathcal{M}$$
We have by construction that $\Gamma_h$ is an invariant curve and $H$ is holomorphic on $\mathcal{M}$. So the Hamiltonian $H$ holomorphic on an open neighbourhood of $\Gamma_h$. We can apply the Morales-Ramis Theorem \ref{thmmorales}, and we get that the variational equation has an Abelian Galois group. We still need to precise on which base field this Galois group is computed; in the Morales-Ramis Theorem, this is the field of meromorphic functions on $\Gamma_h$ (which are meromorphic functions in $q_2(t),p_2(t)$). Using equations \eqref{eqphi1}, \eqref{eqphi2}, we have the following relation for $p_2$
\begin{equation}
p_2^2 =\frac{(q_2^3+1)(hw-\sqrt{2}q_2j)}{wq_2^4} = \frac{(q_2^3+1)(hq_2^2\sin j-8j)}{q_2^6\sin j } \in\mathcal{M}er(\mathcal{S})
\end{equation}
So the field of meromorphic functions on $\Gamma_h$ is just an algebraic extension of degree at most $2$ of $\mathcal{M}er(\mathcal{S})$ and so the identity component of the Galois group of equation \eqref{eq2} will not change if we take $\mathcal{M}er(\mathcal{S})$ instead of this field as base field.

The case of fixed energy is proved by reducing the dynamical system to the hypersurface $H=h$. Indeed, this Hamiltonian field $X_H$ leaves invariant an energy level $H=h$, and so we can define its restriction $X_H\mid_{H=h}$. We lose the Hamiltonian structure and the notion of involution. However, as given by this author in Remark 2 of \cite{87} , if there exists a first integral $I$ on the level $H=h$, then there exists a vector field $X_I=\{I,\cdot\}$ on $H=h$ commuting with the Hamiltonian vector field $X_H\mid_{H=h}$. Thus the field $X_H\mid_{H=h}$ is integrable in the Bogoyavlensky sense \cite{80}: it has one first integral $I$ and two commuting vector fields. We now need to apply a kind of Morales-Ramis Theorem for this $3$-dimensional system $X_H\mid_{H=h}$, and this is given by Theorem 1 of Ayoul-Zung \cite{81}. The normal variational equation of $X_H\mid_{H=h}$ is still equation \eqref{eq2}, and the previous computations are always valid (we already proved that there are no ``exceptional'' energy levels, meaning that the curve $\Gamma_h$ always exists and is never too singular).
\end{proof}

\bigskip

Let us now recall some facts about the monodromy group of a linear differential system. We consider a linear differential system
\begin{equation}\label{eqdiff}
\dot{X}=A(t)X
\end{equation}
with $A$ a matrix whose entries are meromorphic on a Riemann surface $\mathcal{W}$. On $\mathcal{W}$, the equation \eqref{eqdiff} has singularities, corresponding to poles of the entries of $A$. This forms a discrete set $D\subset \mathcal{W}$. We consider a base point $t_0\in\mathcal{W} \setminus D$ and closed curves $\gamma\subset \mathcal{W} \setminus D$ originating in $t_0$. We now compute the resolvant matrix $R(t)$, solution of equation \eqref{eqdiff} with $R(t_0)=I_n$, along the curve $\gamma$. 

Let $R_\gamma$ be the matrix obtained after coming back to $t_0$ for the first time. This matrix is the monodromy matrix of equation \eqref{eqdiff} along $\gamma$. The set of all monodromy matrices forms a group, with the multiplication being the concatenation of closed curves. Moreover, the monodromy group is always contained in the Galois group of equation \eqref{eqdiff} over the base field of meromorphic functions on $\mathcal{W}$.

\bigskip

We have a normal variational equation with transcendental coefficients (and well defined on a Riemann surface). Using a parameter coming from the multivaluation of the coefficients, we will build a ``limit'' equation whose coefficients are rational. We then show that if the Galois group of the normal variational equation is virtually Abelian, then there exists an algebraic relation on the monodromy (Lemma \ref{lemmonderiv}). We prove moreover that this relation goes to the limit, and so is satisfied by the ``limit'' equation.

\bigskip

\begin{lem}\label{lemmonderiv}
Let $\mathcal{M}er(\mathcal{W})$ be the field of meromorphic functions on the Riemann surface $\mathcal{W}$. We consider a differential equation $\dot{X}=A(t) X$ with $A\in M_2(\mathcal{M}er(\mathcal{W}))$. If the Galois group of this system over the base field $\mathcal{M}er(\mathcal{W})$ is virtually abelian, then its monodromy group $G_1$ is such that
$$\forall \; g\in \mathcal{D}^{(2)}(G_1),\;\;g^{60}=id$$
where $\mathcal{D}^{(2)}(G_1)$ is the second derived subgroup of $G_1$.
\end{lem}

\begin{proof}
We know that the monodromy group is a subgroup of the Galois group over the base field of meromorphic functions on $\mathcal{W}$. In dimension $2$, the possible Galois groups (after quotienting them by their center, see \cite{64} Theorem 4.29) are subgroups of the triangular group 
$$\mathbb{B}=\left\lbrace \left(\begin{array}{cc} a&b\\ 0&a^{-1} \end{array}\right) \;\;,a\in \mathbb{C}^*,\;b\in \mathbb{C}\right\rbrace,$$
the infinite diedral group
$$D_\infty=\left\lbrace \left(\begin{array}{cc} a&0\\ 0&a^{-1} \end{array}\right),\left(\begin{array}{cc} 0&a\\ -a^{-1}&0 \end{array}\right) \;\;,a\in \mathbb{C}^*\right\rbrace,$$
or three finite primitiv groups (i.e. $G/\mathcal{Z}(G)\in \{A_4,S_4,A_5\}$). Let us take one of these possible groups and derive it two times (the derivation of a group is the group generated by its commutators). For triangular groups, $D_\infty$, $D_n$, the second derivative produce the identity group
\begin{equation*}\begin{split}
\mathcal{D}(\mathbb{B})=(\mathbb{C},+)\quad \mathcal{D}^{(2)}(\mathbb{B})=id \qquad 
 \mathcal{D}(D_\infty)=(\mathbb{C}^*,\times)\quad \mathcal{D}^{(2)}(D_\infty)=id
 \end{split}\end{equation*}
For the finite primitive groups, we obtain
$$\mathcal{D}^{(2)}(A_4)=id\quad \mathcal{D}^{(2)}(S_4)=\mathbb{Z}_2^2\quad \mathcal{D}^{(2)}(A_5)=A_5.$$
So, in all cases, the elements of the second derivative of the Galois group have always an order dividing $60$. We conclude using the fact that the monodromy group $G_1$ is always a subgroup of the Galois group.
\end{proof}

\begin{lem}\label{thm4}
The normal variational equation \eqref{eq2} does not have a virtually Abelian Galois group over the base field $\mathcal{M}er(\mathcal{S})$ (for any fixed energy $h$).
\end{lem}

\begin{proof}
Assume equation \eqref{eq2} has a virtually Abelian Galois group over the base field $\mathcal{M}er(\mathcal{S})$ (for some fixed energy $h$). We consider the Deck transformation
$$\sigma: \mathcal{M}er(\mathcal{S})\longrightarrow \mathcal{M}er(\mathcal{S}) \quad \sigma\left(j\right)=j+2\pi$$
We now apply $\sigma^l$ on the normal variational equation \eqref{eq2}. This produces the equation
\begin{small}\begin{equation}\label{maineq3}\begin{split}
16\left(h-q_2\left(\delta +\frac{16l\pi}{q_2^3\sin j}\right)\right)(q_2^3+1)(q_2^3-8)^2q_2^4X''+\\
4q_2^3\left(q_2(q_2^6+5q_2^3+112)\left(\delta +\frac{16l\pi}{q_2^3\sin j}\right)+ 24q_2(q_2^3+1)+12h(q_2^3-8)\right)\\
(q_2^3-8) X'-(q_2^3+1)^2\left(q_2^3(q_2^3-32)\left(\delta +\frac{16l\pi}{q_2^3\sin j}\right)-8q_2^3-128\right)X=0
\end{split}\end{equation}\end{small}
We consider the polynomial $P=q_2(q_2^3+1)(q_2^3-8)$, $\nu>0,\epsilon>0$ two real numbers and the compact
$$\mathcal{C}_{\nu,\epsilon}=\left(\left(D(0,\nu)\setminus \left( \cup_{q_2\in P^{-1}(0)} D(q_2,\epsilon)\right)\right)\times D(0,\nu) \right)\cap \mathcal{S}$$
where $D(q_2,\epsilon)$ is a disk with center $q_2$ of radius $\epsilon$. For $l>l_0$ large enough, the singularities of equation \eqref{eq2} are not in $\mathcal{C}_{\nu,\epsilon}$, because for large $l$, the singularities of equation \eqref{eq2} inside $D(0,\nu)^2 \cap \mathcal{S}$ are converging to the roots of $P$ (and singular points of $\mathcal{S}$ which are also roots of $P$). Now let consider four paths $\gamma_1,\dots,\gamma_4$ on $\mathcal{S}$ outside the roots of $P$ and the commutator
\begin{equation}\label{commutator}
[[R_{\gamma_1},R_{\gamma_2}],[R_{\gamma_3},R_{\gamma_4}]]^{60}
\end{equation}
where $R_\gamma$ is the resolvant matrix of equation \eqref{maineq3} along $\gamma$. As equation \eqref{eq2} has a virtually Abelian Galois group over the base field $\mathcal{M}er(\mathcal{S})$, so is the case of equation \eqref{maineq3}. Using Lemma \ref{lemmonderiv}, the monodromy group $G_1$ of equation \eqref{maineq3} is such that
$$\forall \; g\in \mathcal{D}^{(2)}(G_1),\;\;g^{60}=id$$
Thus the commutator \eqref{commutator} is equal to identity for any $l$.

There exist $\nu,\epsilon$ such that $\gamma_1,\dots,\gamma_4\subset \mathcal{C}_{\nu,\epsilon}$. We divide \eqref{maineq3} by its dominant term and we take the limit $l\longrightarrow \infty$. The equation is then converging to
\begin{equation}\label{maineq2}\begin{split}
16q_2^2(q_2^3+1)(q_2^3-8)^2y''-4q_2(q_2^6+5q_2^3+112)(q_2^3-8)y'+\\
(q_2^3+1)^2(q_2^3-32)y=0
\end{split}\end{equation}
The resolvant matrix $R$ is smooth when $l\longrightarrow\infty$ on the compact set $\mathcal{C}_{\nu,\epsilon}$. So, for $l$ large enough, the resolvant matrices $R$ are smooth along $\gamma_i$, and so they are converging to monodromy matrices along curves $\gamma_i$ of equation \eqref{maineq2}. So is the commutator \eqref{commutator}, which is equal to $id$ by hypothesis. So the monodromy group $G_2$ of equation \eqref{maineq2} is such that
$$\forall \; g\in \mathcal{D}^{(2)}(G_2),\;\;g^{60}=id$$
Analysis of singularities of equation \eqref{maineq2} shows that all its singularities are regular, and then that the equation is Fuchsian. This implies that the Galois group $G_3$ of equation \eqref{maineq2} over the differential field $K_0$ of rational functions on $\mathcal{S}$ is exactly the Zariski closure of the monodromy group $G_2$. Thus we obtain also the property
$$\forall \; g\in \mathcal{D}^{(2)}(G_3),\;\;g^{60}=id$$
This implies that the identity component of $G_3$ is solvable. We have
$$\hbox{Gal}_{\hbox{diff}}(K_0/\mathbb{C}(t))=D_\infty$$
So the Galois group over $\mathbb{C}(t)$ of equation \eqref{maineq2} is a solvable extension of the Galois group $G_3$. Thus the Galois group over $\mathbb{C}(t)$ of equation \eqref{maineq2} should have a solvable identity component. Using Kovacic algorithm on equation \eqref{maineq2}, we prove that its Galois group over $\mathbb{C}(t)$ is $SL_2(\mathbb{C})$, which is connected and not solvable.

So equation \eqref{eq2} has not a virtually Abelian Galois group over the base field $\mathcal{M}er(\mathcal{S})$ for any fixed energy $h$.
\end{proof}

\noindent
Using Theorem \ref{thm3} and Lemma \ref{thm4}, the Hamiltonian $H$ has not an additional first integral in $\mathcal{M}er(\mathcal{M})$, even if restricted to a single energy level $H=h$. This implies the main Theorem \ref{thm0}.

\section{Conclusion}
It is not so rare that a transcendental function appear in the variational equation. To prove non-integrability, we need to study its Galois group. In general it is possible to avoid it by just taking a particular orbit for which this case does not occur, as done for example in \cite{32}. But we see that in fact it probably produces even stronger integrability conditions. Such a study is not so much more difficult when through a limiting process, such an equation induces a ``limit'' equation with coefficients in $\mathbb{C}(t)$. The integrability condition on this ``limit'' equation is that the Galois group should be virtually solvable. The original condition was virtual abelianity, but this does not change anything in practice (aside in higher variational equations). And moreover, it is only a necessary criterion, if it was met, we could produce other conditions by making an asymptotic expansion in the ``multivaluation parameter'' (which corresponds here to apply the sheave translation $\sigma$). In the case where these transcendental extensions are not avoidable, this approach completely make sense because the equivalent of Kovacic algorithm for such equations is not implemented yet.

In the same problem of Riemann ellipsoid motion, zero angular momentum is only one case. A (probably) complete list of integrable cases is given in \cite{28}, and it could maybe be possible to prove the non-integrability of the other cases using this approach.

\label{}

\bibliographystyle{plain}
\bibliography{bibthese}

\end{document}